\documentclass[10pt]{article}
\usepackage{ametsoc2col}
\newcommand{\myabstract}{It is of crucial importance to be able to identify the location of atmospheric pollution sources in our planet. Global models of atmospheric transport in combination with diverse Earth observing systems are a natural choice to achieve this goal. It is shown that the ability to successfully reconstruct the location and magnitude of an instantaneous source in global chemical transport models (CTMs) decreases rapidly as a function of the time interval between the pollution release and the observation time. A simple way to quantitatively characterize this phenomenon is proposed based on the effective --undesired-- numerical diffusion present in current Eulerian CTMs and verified using idealized numerical experiments. The approach presented consists of using the adjoint-based optimization method in a state-of-the-art CTM, GEOS-Chem, to reconstruct the location and magnitude of a realistic pollution plume for multiple time scales. The findings obtained from these numerical experiments suggest a time scale of 2 days after which the accuracy of the adjoint-based optimization methodology is compromised considerably in current global CTMs.  In conjunction with the mean atmospheric velocity, the aforementioned time scale leads to an estimate of a length scale of about $1700km$, downwind from the source, beyond which measurements, in conjunction with current global CTMs, may not be successfully utilized to reconstruct continuous-in-time sources.  The approach presented here can be utilized to characterize the capabilities and limitations of adjoint-based optimization inversions in other regional and global Eulerian CTMs. }
\begin{document}
%
%
\title{\textbf{\large{Quantifying the loss of information in source attribution problems using the adjoint method in global models of atmospheric chemical transport. }}}
%
%
\author{\textsc{Mauricio Santillana,}
				\thanks{\textit{Corresponding author address:} 
				Mauricio Santillana, 29 Oxford St., Cambridge, MA 02138. 
				\newline{E-mail: msantill@fas.harvard.edu}}\\
\textit{\footnotesize{Harvard University School of Engineering and Applied Sciences, Cambridge, Massachusetts}}
}
%
\ifthenelse{\boolean{dc}}
{
\twocolumn[
\begin{@twocolumnfalse}
\amstitle

\begin{center}
\begin{minipage}{13.0cm}
\begin{abstract}
	\myabstract
	\newline
	\begin{center}
		\rule{38mm}{0.2mm}
	\end{center}
\end{abstract}
\end{minipage}
\end{center}
\end{@twocolumnfalse}
]
}
{
\amstitle
\begin{abstract}
\myabstract
\end{abstract}
\newpage
}
\section{Introduction}
Understanding and quantifying the fate of anthropogenic and natural emissions of chemically active gases in the atmosphere is an important endeavour in our changing climate. In particular, many efforts have been directed to the quantification of large-scale spatial and temporal variations of sources and sinks of gases such as carbon dioxide, methane and nitrous oxide, due to their direct role in changing the radiative properties of the atmosphere, see \citet{bib:Kas00, bib:Ent02, bib:Cia10} and the multiple references therein. Having reliable estimates of these sources and sinks is of crucial relevance in the evaluation of global policies such as the Kyoto protocol, designed to curb emissions of green-house gases. \\

Global chemical transport models (CTMs) have been designed to simulate the dynamics of the concentration fields of chemicals under the influence of atmospheric transport and chemical reactions. In conjunction with Earth observing systems such as satellite retrievals and monitoring stations, CTMs provide a natural modeling framework to estimate the strength and location of such chemical sources in a top-down fashion. In this context, CTMs are expected to have the capability of identifying pollution sources within the time scales of relevance to global atmospheric chemical transport events: weeks, for intercontinental transport of pollutants, and months to a year, for inter-hemispheric transport. It is shown in this study that modeling frameworks that use the adjoint-based optimization method to identify the strength and location of pollution sources, as implemented in current CTMs, may not achieve their goal in relevant time scales. \\

Our ability to reconstruct a pollution source using global CTMs depends mainly on three factors: ({\it i}) the quality of our measurement systems, ({\it ii}) the appropriateness of our mathematical model 
and our inversion approach, and ({\it iii}) our ability to numerically approximate the mathematical model and the inversion approach using computers \citep{bib:Ent02}. In this study, it is assumed that a perfect observation system is in place (no noise) at every location in the atmosphere, thus ignoring errors coming from ({\it i}), and it is also assumed that the mathematical model (equations (\ref{eq_transp_chem})) approximates well the dynamics of the atmosphere, thus ignoring the issue of the appropriateness of the mathematical model in ({\it ii}). Specifically,  the efficacy of the adjoint-based optimization methodology is studied in a practical and realistic computational framework ({\it iii}). \\

Eulerian global CTMs simulate the chemical composition of the atmosphere by numerically solving a set of coupled partial differential equations of the form,
\begin{equation}
\frac{\partial C_i}{\partial t}+\mathbf{u}\cdot\nabla(C_i)=R_i+s_i(\mathbf{x}, t)
\label{eq_transp_chem}
\end{equation}
where $C_i$ is the concentration of chemical species $i$, $\mathbf{u}$ is the wind velocity field (obtained from a global atmospheric circulation model), $R_i$ is the effective chemical production rate (typically a function of the mass fractions 
of other chemicals), and $s_i$ describes local emissions and non-chemical sinks. \\

A source attribution problem in a global CTM can be stated as follows, given a set of observations of the state of the atmosphere (typically concentration fields) and atmospheric wind fields for a time period $T=t_f-t_0>0$, find the best emission and deposition rates, $s_i(\mathbf{x}, t)$, for the time period $T$, such that when they are utilized as parameters in the CTM, the simulated chemical concentrations $C_i(t)$ are consistent with observations. Finding the solution  (the best emission and depostion rates) to this inverse problem is achieved by minimizing the misfit between observations and model results, frequently using gradient-based iterative optimization methods. The adjoint-based optimization method studied here, is a widely accepted methodology to calculate the gradients needed to identify the direction of steepest descent in the minimization of a cost function that typically quantifies the misfit between model simulations and measurements. This methodology has been widely utilized in the reconstruction of emission sources using global CTMs, see for example \citep{bib:Hen07, bib:Zha09, bib:Hen09, bib:Kop09, bib:Kop10, bib:Zha11, bib:Wec11}.\\

In practical terms, the calculation of the gradient at each iteration using the adjoint is achieved by integrating the CTM forward-in-time and back-in-time with reversed winds \citep{bib:Talagrand87}. Thus, the efficacy of the adjoint-based optimization methodology depends directly on the properties of the numerical schemes utilized in the CTM, along with the properties of atmospheric flow. In particular, the presence of numerical errors in Eulerian advection schemes utilized to integrate equations (\ref{eq_transp_chem}) cause the true solution of equations (\ref{eq_transp_chem}) and the (computational) numerical approximation to the solution, to diverge as the integration time $t$ increases \citep{bib:Lev02}. This numerical divergence also takes place between the true adjoint and the numerical calculation of the adjoint, and thus, in the approximation of the gradient to be utilized in the minimization process. As a consequence, it is expected that the efficacy of the numerical solution of source attribution problems, using the adjoint-based optimization method, will lead to unreliable source reconstructions for large enough simulation times. \\

It is shown in this work that unreliable source reconstructions may happen in global Eulerian CTM simulation times of about 2 days, a considerably short time scale in the context of global atmospheric dynamics (compare it to a week, a typical time scale for inter-continental transport of pollution plumes). An idealized accident-type source reconstruction problem is utilized to investigate this and in general the capabilities and limitations of the adjoint-based optimization methodology. The aforementioned time scale in combination with the mean atmospheric velocity leads to a length scale of about $1700km$ (or about 3 grid boxes downwind), downwind from the source, beyond which measurements may not be successfully utilized to solve source attribution problems in current global CTMs. A simple way to estimate this crucial time scale is proposed here for any regional or global Eulerian CTM, based on the ``effective'' numerical diffusion present in the model, which happens to be much greater than the numerical diffusion determined by the order of the advection scheme, due to the chaotic properties of atmospheric flow.\\


This paper is organized as follows, the idealized source attribution problem is introduced in section \ref{sec:source}, the adjoint-based optimization method, as well as the technical details of its numerical implementation are presented in section \ref{sec:adjoint}. The numerical results of our computational experiments along with their implications in real-world applications, including data assimilation methodologies, are discussed in section \ref{sec:numerical_results}.  



\section{A source attribution problem}
\label{sec:source}

An idealized accident-type release source attribution problem for an inert gas is studied here in order to evaluate the efficacy of the adjoint-based methodology to reconstruct local-in-time and space perturbations. Thus, from here on $R_i=0$ in equations (\ref{eq_transp_chem}). An observation system capable of providing exact (synthetic) observations of the atmosphere everywhere at time $t_f$ (obtained directly from a forward simulations using the CTM) is assumed. The problem consists of finding the best field of instantaneous emissions at time $t=t_0$, $s_i(\mathbf{x}, t_0)$, such that the misfit between observations and the simulation results, at time $t_f$, is minimized. The  problem is simplified by placing the instantaneous release at a vertical height of approximately $4$ km above the ground at time $t=0$. This choice minimizes the effects of planetary boundary layer mixing processes, as represented (by sub-grid parametrizations) in global CTM \citep{bib:Lin08}. By having the same number of observations, $p$ (all grid boxes in the domain at $t=t_f$), as the number of discrete instantaneous sources to be determined, $N$ (all grid boxes in the domain at $t=0$), the usual ill-posed nature of the inversion problem caused by the under-determination of parameters when $p\ll N$ in real-life Earth observations systems, is removed \citep{bib:Boc05a, bib:Boc05b}. This fact, in combination with the numerical stability provided by the use of the numerical approximation of the continuous adjoint (see section \ref{sec:adjoint}) as a means to calculate the gradients in the minimization process, makes it possible for the idealized experiments under examination to be solved without the need of prior information about the plume or a regularization technique \citep{bib:Boc05a, bib:Boc05b}. \\

\begin{figure}[t]
\noindent\includegraphics[width=.5\textwidth]{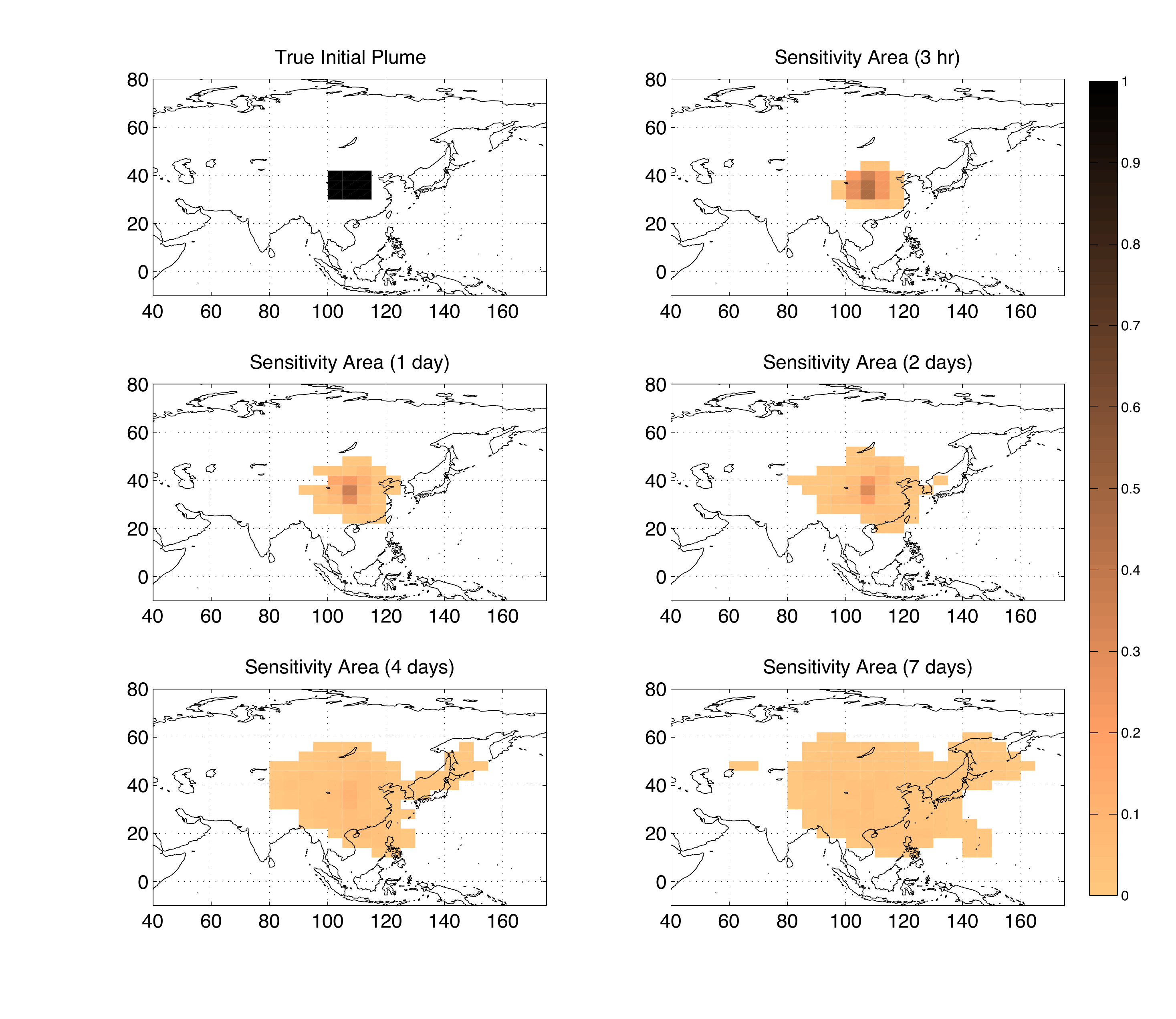}
\caption{Area of influence at $\sim 4$ km of altitude (20th pressure level) as a function of simulation time, as reconstructed by the first iteration of the adjoint-based optimization method for $4^{\circ}\times 5^{\circ}$ numerical simulations. Latitude vs Longitude. }
\label{fig_sensitivity4x5}
\end{figure}

\begin{figure}[t]
\noindent\includegraphics[width=.5\textwidth]{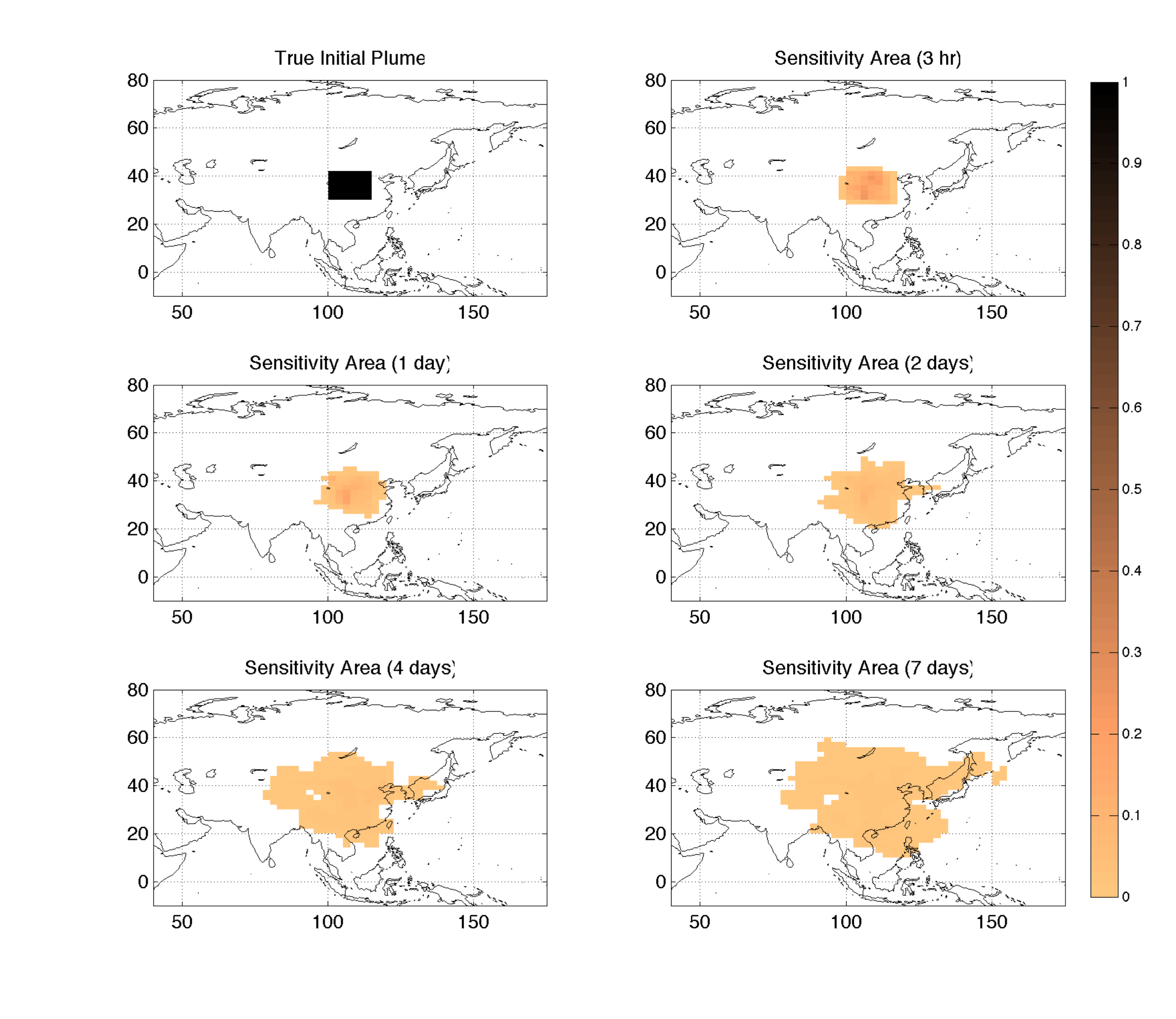}
\caption{Area of influence at $\sim 4$ km of altitude (20th pressure level) as a function of simulation time, as reconstructed by the first iteration of the adjoint-based optimization method for $2^{\circ}\times 2.5^{\circ}$ numerical simulations. Latitude vs Longitude. }
\label{fig_sensitivity2x25}
\end{figure}

The aforementioned problem is solved numerically, using a three dimensional state-of-the-art CTM, GEOS-Chem \citep{bib:Bey01}. GEOS-Chem is a state-of-the-art CTM, driven by GEOS-5 analyzed meteorological data from the NASA Global Modeling and Assimilation Office (GMAO). The GEOS-Chem adjoint-based optimization was developed by \citep{bib:Hen07}.  In order to produce synthetic observations, an inert tracer plume is propagated for different simulation times in GEOS-Chem, with $4^{\circ}$ x $5^{\circ}$ and $2^{\circ}$ x $2.5^{\circ}$ horizontal resolutions, using the high-order advection scheme native of GEOS-Chem \citep{bib:Lin96}. The true initial plume concentration was set to be twice that of the constant global background in all other grid cells, with a horizontal extent of approximately $1330$ km $\times$ $1680$ km (3 x 3 grid cells), and a vertical extent of one vertical pressure level at 4km (20th pressure level) above the ground. The horizontal location of the initial plume is shown in Figures \ref{fig_sensitivity4x5} and \ref{fig_sensitivity2x25} (Top left plot). Through successive iterations, the adjoint-based optimization methodology  generates an ``optimized'' initial condition estimate. The maximum number of iterations was set to 99 to keep our experiments within the practical limit set by our computational resources. We performed six numerical experiments for each of the two spatial resolutions, assimilating data after 3, 12, 24, 48, 96, and 168 hours.\\

Since any continuous emission field in time can be thought of as a sum of instantaneous releases, the structural results of our numerical investigation are valid in constant emission flux reconstructions and more generally, in data assimilation frameworks that use Eulerian CTMs and the adjoint-based optimization method to estimate source locations and/or magnitudes. \\


The idealized numerical problem studied here is called an ``identical twin'' experiment in the inverse modeling community \citep{bib:Kai04}, since synthetic observations are produced using the same forward model utilized in the inversion approach. It is well known that assessing the performance of inverse modeling strategies using an identical twin experiment leads to optimistic results and constitutes what the inverse modeling community calls an ``inverse crime'' \citep{bib:Kai04}. Indeed, in identical twin experiments, the exact solution of the minimization problem lives in the parameter space (the space of initial conditions in our specific problem) of the modeling framework, and thus, with an appropriate methodology one should be able to recover the optimal solution. From this perspective, the idealized numerical experiments studied here assess only the algorithmical capabilities and limitations of the adjoint-based inversion methodology, as implemented in CTMs, and can be interpreted as optimistic since the further complications of real-life source reconstructions are ignored.

\section{The adjoint-based optimization method}
\label{sec:adjoint}
Inverse modeling and data assimilation approaches involve finding an optimal set of model parameters that best match observations (\cite{ bib:Ent02, bib:Ben02, bib:Wun06}). In order to find this optimal set of parameters (in our source attribution problem we find the best initial condition), one minimizes the distance or misfit between the model results and observations. This minimization can be performed using a gradient based approach. As shown in \citep{bib:Talagrand87} the ``adjoint equations''  of the model can be used to compute explicitly the functional derivative (or the gradient) of the distance function, between observations and model results, with respect to the initial conditions.  The computation of one gradient requires one forward-in-time integration of the full model equations over the time interval on which the observations are available, followed by one backward-in-time integration of the adjoint equations. Subsequent gradients required for a steepest descent algorithm in a minimization routine are computed in a similar fashion.\\  

When $R_i=0$ in equations (\ref{eq_transp_chem}), the forward model becomes a linear transport problem with sources and sinks.  The adjoint operator of the transport operator is again a linear transport operator with reversed winds and reversed time integration \citep{bib:Boc05b, bib:Hou06}. A subtle point arises when implementing the adjoint-based optimization method numerically, since one faces the dilemma of whether using the adjoint of the discrete approximation to the transport operator or the discretization of the continuous adjoint \citep{bib:Vuk01, bib:Boc05a, bib:Boc05b, bib:Hou06, bib:Liu08}. This issue has been previously studied not only in the context of atmospheric chemistry modeling \citep{bib:Gou11}, but in the context of oceanography and weather modeling \citep{bib:Sir97}. While the adjoint of the discrete approximation of the transport may produce better point-wise results in sensitivity analyses, it can be very unstable numerically and potentially lead to non-physical sensitivity results, specially when using stabilized high-order advection schemes that are not time symmetric \citep{bib:Hou06}.
The superior numerical stability and global convergence properties of the numerical approximation to the continuous adjoint make it a preferred choice in practical implementations of atmospheric transport. As a consequence, the adjoint operator is calculated using the same transport routine as in the forward simulations \citep{bib:Vuk01, bib:Hou06, bib:Hen07, bib:Gou11} and thus, the efficacy of this inversion approach inherits the capabilities and limitations of the transport routine.\\

\subsection{Numerical diffusion versus atmospheric eddy diffusion}

Numerical errors in transport routines of Eulerian global CTMs have been shown to be significant, compromising the quality of simulation results, at spatial resolutions larger than $1^{\circ} \times1^{\circ}$ even with high-order advection schemes (Wild and Prather, 2006). Moreover, the decay and spatial broadening observed in Eulerian global CTM simulations of realistic pollution plumes signal the presence of excessive numerical diffusion. This undesired numerical diffusion is shown to be much higher than expected by the order of the advection scheme, and shown to be nearly governed by the theoretical upper limit,
set by the value of the local finite-time Lyapunov exponent, almost independently of the grid-box size, in \citet{bib:Ras11} for currently practical spatial resolutions. As a consequence, high values of numerical diffusion are observed in regions with high values of finite-time Lyapunov exponents, \textit{i.e.} in regions where wind shear is high and where chaotic mixing (and thus where loss of information) is pronounced. 
This excessive numerical diffusion will eventually compromise our ability to recover the plume information (location and magnitude) in the multiple iterations of the adjoint-based optimization process. \\

Some amount of diffusion is to be expected in any atmospheric transport simulation. Indeed, turbulent mixing, consequence of the multi-scale chaotic velocity field in the atmosphere, can be thought of and is often parametrized as an eddy diffusion term in meso-scale Eulerian atmospheric simulations and represented by stochastic perturbations on the velocity field in Lagrangian global CTMs \citep{bib:Gif82, bib:Sil90, bib:Mau98, bib:Pis09}. Studies combining Lagrangian trajectories and aircraft measurements of actual passive tracer plumes have placed the values of this eddy diffusivity around $10^4$ m$^2$ s$^{-1}$ in the horizontal, and $0.3 - 1$ m$^2$ s$^{-1}$ in the vertical \citep{bib:Pis09}. In Eulerian global CTMs, however, the values of the ``effective'' horizontal numerical diffusion are estimated to be approximately $10^5$ m$^2$ s$^{-1}$ (Pisso et al., 2009), an order of magnitude larger than those estimated by aircraft studies. Specifically, this horizontal numerical diffusion in GEOS-Chem can be estimated\footnote{For this estimate we assumed that the decay rate for a spatial resolution of $4^{\circ}\times 5^{\circ}$ is governed by the Lyapunov exponent $\lambda$, thus $\alpha\approx 9^{-5}$ sec$^{-1}$, the plume width $W\sim 10^6$ m, and the characteristic length scale where the concentration of the plume decays to zero $r_b\approx \Delta x= 4\times 10^5$ m.} to be about $3.6\times 10^6$ m$^2$ s$^{-1}$, using equation (13) in (Rastigeyev et al., 2010) for realistic plumes in the mid-latitudes of sizes of $\sim 1000$ km. This value is in fact two orders of magnitude larger than those reported in aircraft studies which is why horizontal eddy diffusion terms are not considered, as sub-grid parametrizations, in Eulerian global CTMs. Studying the effect of this excessive numerical diffusion on source attribution problems is the main motivation of this study. \\



The effective (numerical) equations governing the dynamics of a pollution plume in global Eulerian CTMs can be incorporated in the mathematical model (\ref{eq_transp_chem}) as:
\begin{equation}
\frac{\partial C_i}{\partial t}+\mathbf{u}\cdot\nabla(C_i)=\nabla\cdot(\mathbf{D_h}\;\nabla(C_i))+R_i+s_i(\mathbf{x}, t)
\label{eq_transp_chem_diff}
\end{equation}
where the ``effective'' numerical diffusivity matrix $\mathbf{D_h}$ is, in fact, a function of the grid size, $\Delta \mathbf{x}$, the time step, $\Delta t$, the characteristics of the flow field given by $\mathbf{u}(\mathbf{x}, t)$, and the spatial extent of the pollution plume with respect to the grid size \citep{bib:Lev02, bib:Maj99, bib:Ras11}. 

\section{Numerical Results}
\label{sec:numerical_results}

A gradual process to investigate the impact of the effective numerical diffusion in the solution of the aforementioned idealized source attribution problems was designed. First,  one iteration of the adjoint-based optimization was performed in order to identify the geographic area of influence, {\it i.e.} where information came from, as calculated by the numerical advection scheme in our numerical experiments for multiple time-scales. Since no prior knowledge of the pollution plume was assumed, the first meaningful iteration of the adjoint-based optimization consists of utilising the resulting concentration field of the forwardly advected idealized plume as an initial condition, and integrate it back-in-time using the advection scheme, with reversed winds, for exactly the same time period $t_f-t_0=t_f$.  In the absence of numerical errors, and due to the linearity of the transport operator in equations (\ref{eq_transp_chem}),
this procedure should reconstruct the location and magnitude of the original plume at time $t=0$. 

\subsection{Geographic region of influence and loss of information}
\label{sec:geo_inf}

The obtained results from this single iteration of the adjoint-based optimization, displayed in Figures \ref{fig_sensitivity4x5} and \ref{fig_sensitivity2x25} for spatial resolutions $4^{\circ}\times 5^{\circ}$ and $2^{\circ}\times 2.5^{\circ}$ respectively, 
show the numerical reconstruction of the geographic region of influence, $A_h(t_0)$, also called the numerical sensitivity region, at $\sim 4$ km altitude (at time $t=0$). Note that for a simulation time $t_f=3h$, the numerical region of influence coincides with the original region of the plume in both spatial resolutions. As the simulation time $t_f$ (time between release and observations) increases, the numerical sensitivity region broadens more and more (the boundary of this geographic region was drawn where the reconstructed source was $100$ times smaller than the original excess concentration). The broadening --a direct consequence of the numerical diffusion-- is slightly anisotropic reflecting the underlying structure of the atmospheric winds characterized by the jet stream in mid-latitudes. Using estimates of the local finite-time Lyapunov exponent for mid-latitudes in GEOS-Chem from \cite{bib:Ras11}, for $4^{\circ}\times 5^{\circ}$ and $2^{\circ}\times 2.5^{\circ}$ spatial resolutions, we can infer approximate values of $D_h\sim 3.5 \times 10^{6}$ m$^2$ s$^{-1}$  and $D_h\sim 1.4 \times 10^{6}$ m$^2$ s$^{-1}$ for the numerical horizontal diffusion for each resolution respectively, in the chosen geographic location of the idealized plume.\\

Assuming horizontal isotropic diffusion only, an estimate of the ``diffusive deformation'' of the numerical region of influence can be calculated by adding the length scale $\sqrt{D_h t}$ to each the boundary lengths of the plume at time $t$. We can thus expect that the area of the numerical region of influence will have increased its size nearly four-fold in $4^{\circ}\times 5^{\circ}$ simulations and three-fold in $2^{\circ}\times 2.5^{\circ}$, after a day of forward simulation (and a day of back integration),  since the tracer concentration will have spread $\sim 1000$ km on each of the four boundaries of the plume.
As shown in the Sensitivity Area (1 day) panel of Figures \ref{fig_sensitivity4x5} and \ref{fig_sensitivity2x25}, this estimate captures well the behaviour of the numerical results. For two-days and four-days forward simulation times, this estimate suggests a 4.5-fold and 6-fold (for $4^{\circ}\times 5^{\circ}$ simulations) and 4-fold and 5.4-fold  (for $2^{\circ}\times 2.5^{\circ}$ simulations) broadening of the numerically reconstructed plume, respectively. In fact, the introduced deformation estimates of the area of the numerical region of influence appear to be a lower limit of the actual one shown in our numerical experiments in the panels of Figure \ref{fig_sensitivity4x5}. The disagreement, however, is not significant.  Moreover, it is shown in the next section that this estimates of the --undesired-- numerical spatial deformation of the original plume can be used to estimate the loss of information of the full adjoint-based optimization.\\

While this estimate of the deformation of the numerical sensitivity region, using only (spatially and time averaged) horizontal isotropic broadening, may appear incomplete (local finite-time Lyapunov exponent values, from which we estimated the spatially and time averaged diffusivity constant in the diagonal of $\mathbf{D_h}$, change as a function of time and region of space), it seems appropriate when calculating the effective horizontal broadening of the horizontal component of the boundaries of the plume (and thus the deformed area of the plume) and it is consistent with the findings shown in Figure 3 of \citep{bib:Ras11}, where 2D plume decay estimates (and thus plume broadening estimates) approximate well the plume decay even in 3D simulations in GEOS-Chem. \\

Note that when the reconstructed numerical area of influence differs significantly from the original area of influence, as is the case for $t> 2$ days (Figures \ref{fig_sensitivity4x5} and \ref{fig_sensitivity2x25}), the gradient information provided by the adjoint operator is of poor quality and as a consequence, the adjoint-based optimization fails to produce substantial improvements in subsequent iterations (Figures \ref{fig:cost2x25} and \ref{fig:cost4x5}). This statement is investigated quantitatively in the following section.

\begin{figure}[t]
\noindent\includegraphics[width=.5\textwidth]{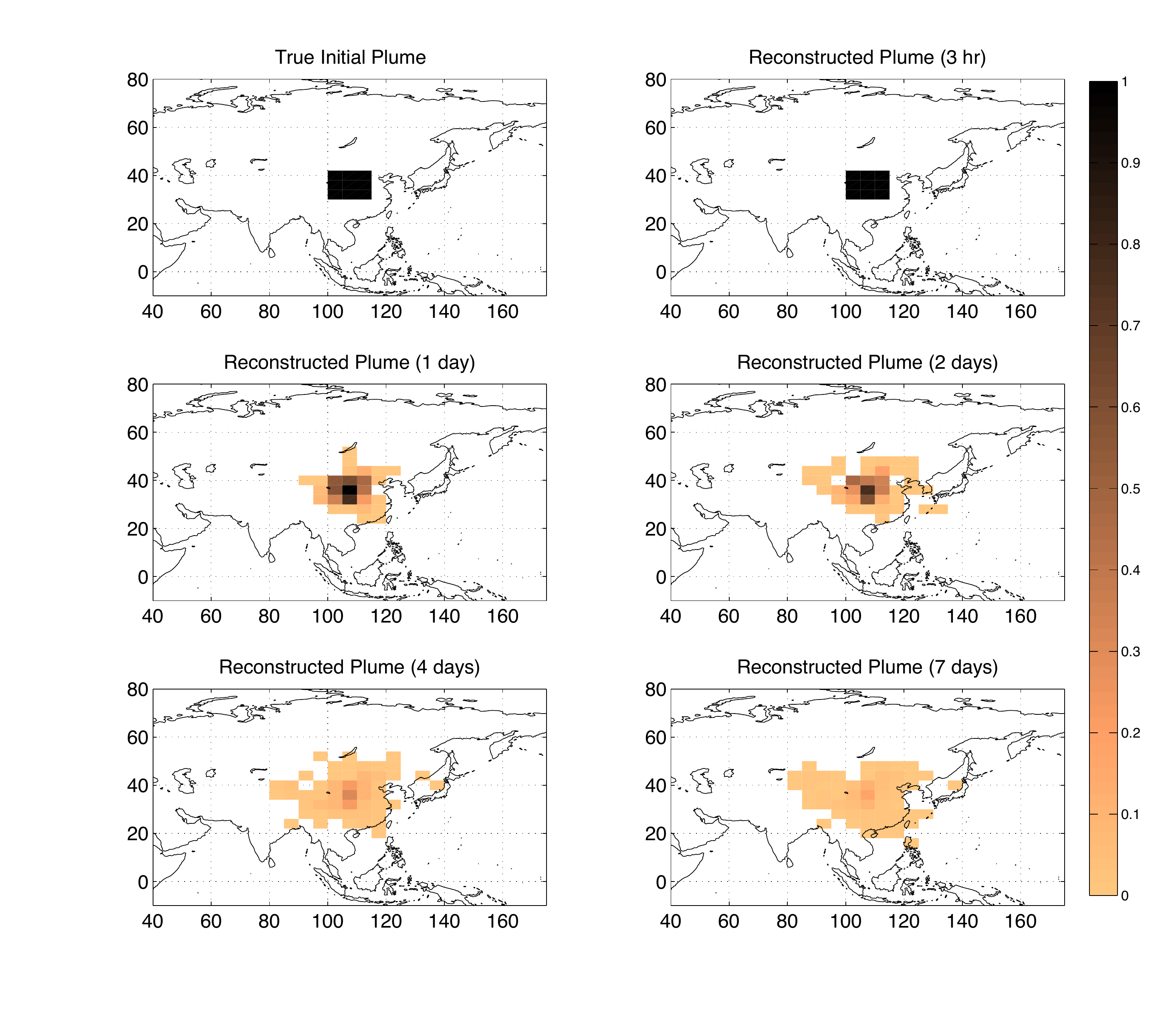}
\caption{Reconstructed plumes at $\sim 4$ km of altitude (20th pressure level) as a function of simulation time, as reconstructed by the adjoint-based optimization method for $4^{\circ}\times 5^{\circ}$ numerical simulations. Latitude vs Longitude.  }
\label{fig_reonstructions_4x5}
\end{figure}

\begin{figure}[t]
\noindent\includegraphics[width=.5\textwidth]{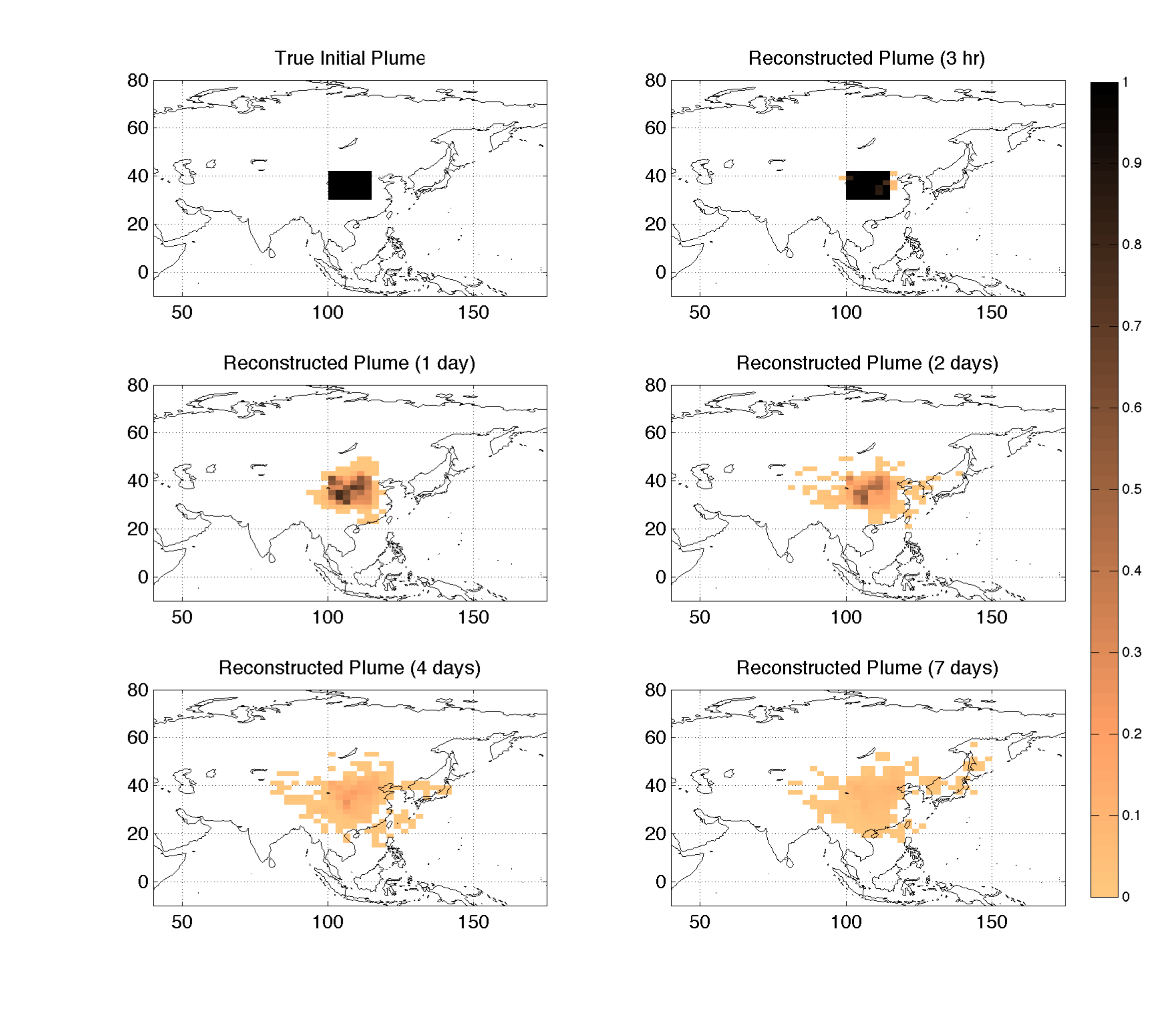}
\caption{Reconstructed plumes at $\sim 4$ km of altitude (20th pressure level) as a function of simulation time, as reconstructed by the adjoint-based optimization method for $2^{\circ}\times 2.5^{\circ}$ numerical simulations Latitude vs Longitude.}
\label{fig_reonstructions_2x25}
\end{figure}

\begin{figure}[t]
\noindent\includegraphics[width=.5\textwidth]{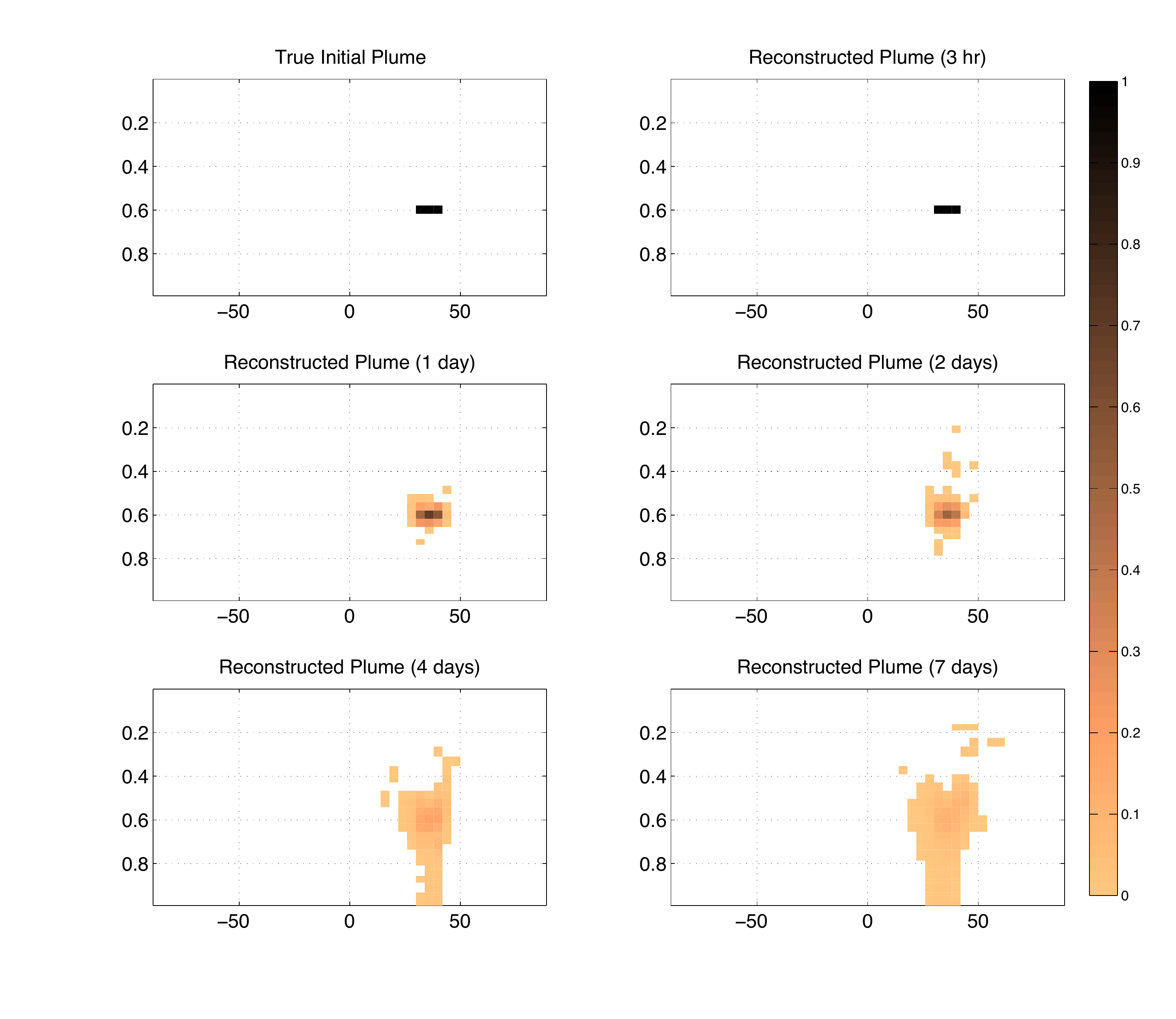}
\caption{Normalized, zonally-integrated reconstructed plumes as a function of simulation time, as reconstructed by the adjoint-based optimization method for $4^{\circ}\times 5^{\circ}$ numerical simulations. Latitude vs Pressure level (atm)}
\label{fig_reconstr_4x5_lat_vs_alt}
\end{figure}

\begin{figure}[t]
\noindent\includegraphics[width=.5\textwidth]{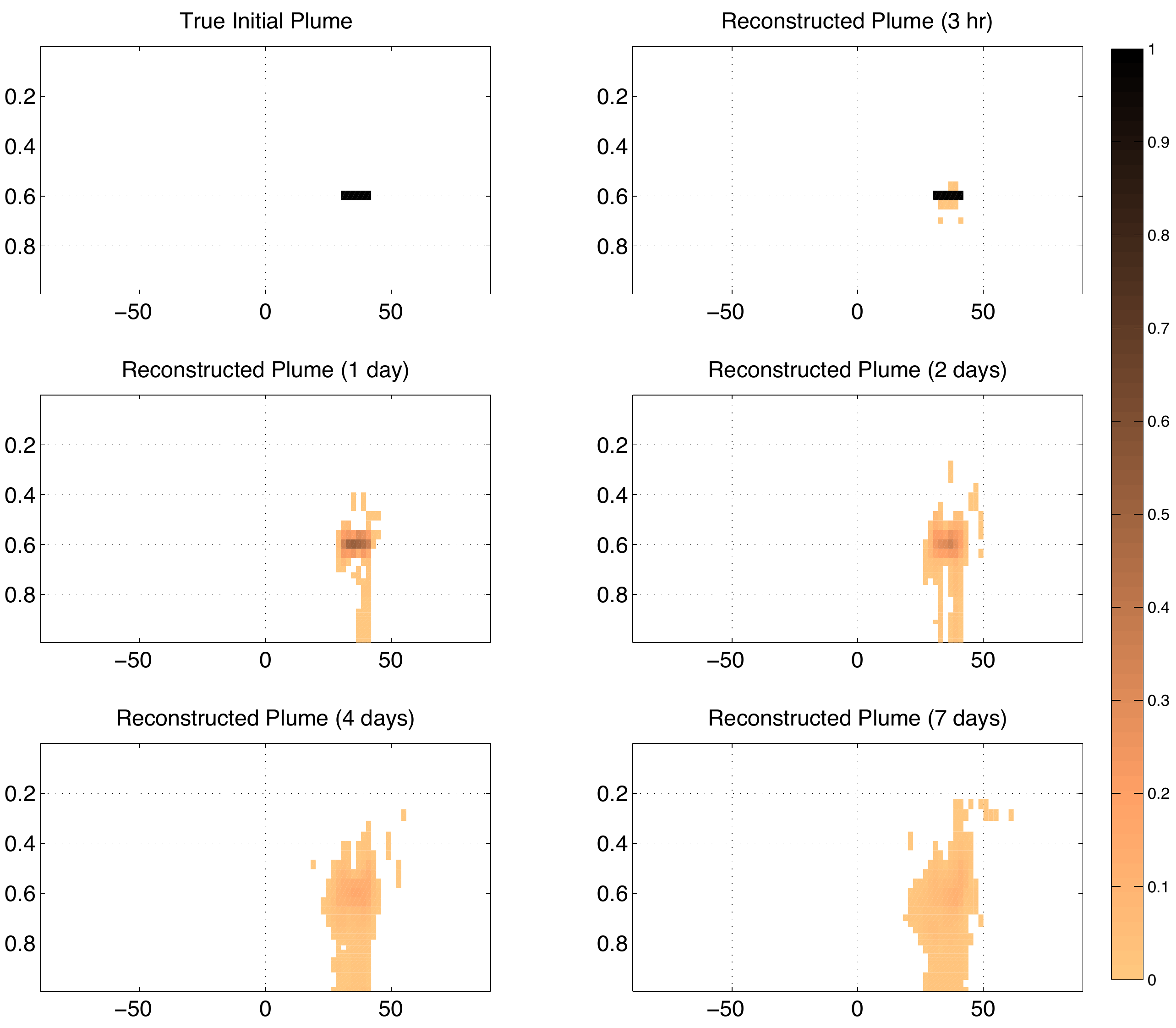}
\caption{Normalized, zonally-integrated reconstructed plumes as a function of simulation time, as reconstructed by the adjoint-based optimization method for $2^{\circ}\times 2.5^{\circ}$ numerical simulations. Latitude vs Pressure level (atm)}
\label{fig_reconstr_2x25_lat_vs_alt}
\end{figure}

\begin{figure}[t]
\noindent\includegraphics[width=.5\textwidth]{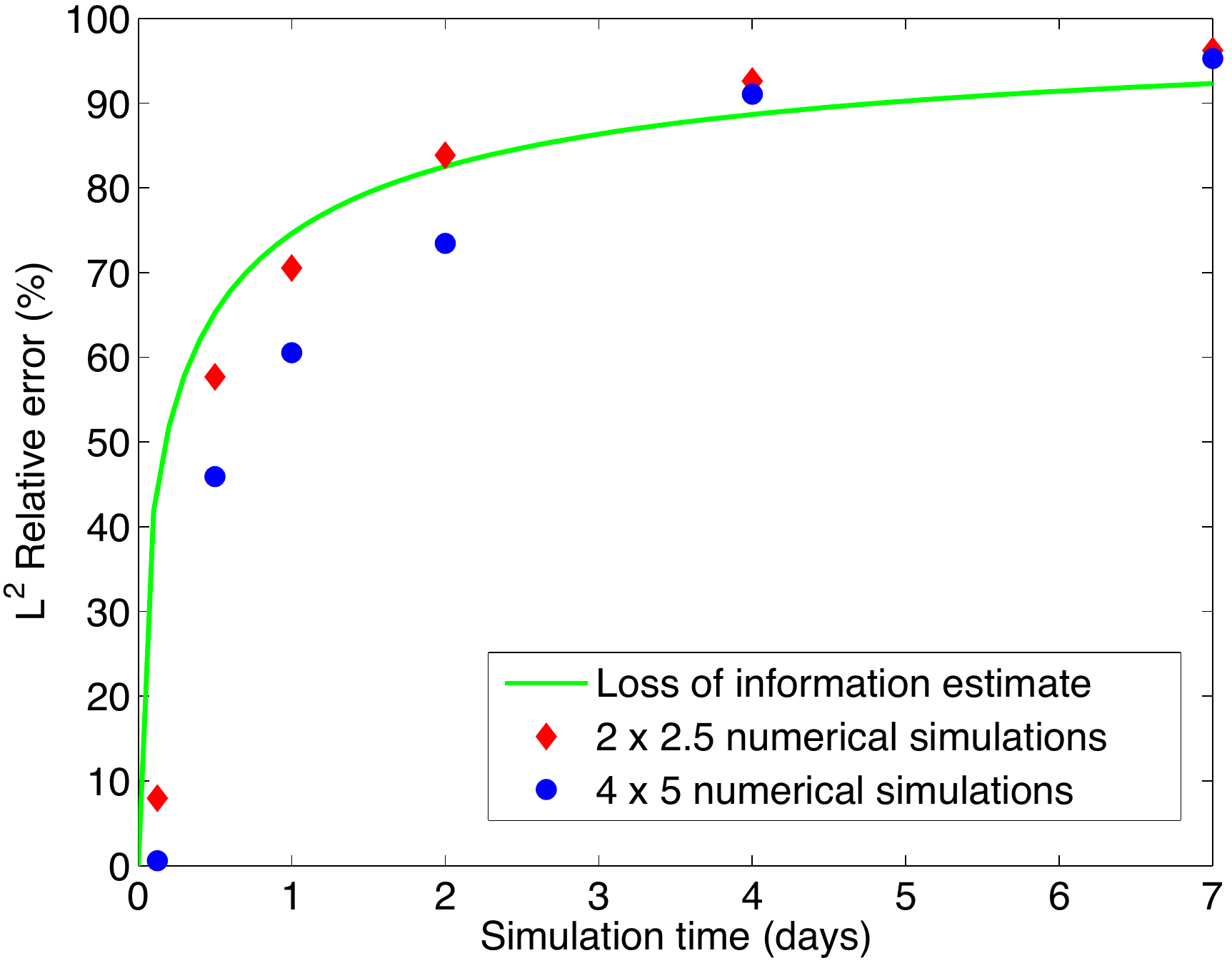}
\caption{Relative $L^2$-error ($\%$) between the reconstructed initial condition and the true solution as a function of the simulation time for both $4^{\circ}\times 5^{\circ}$ and $2^{\circ}\times 2.5^{\circ}$ spatial resolutions. The estimate on the loss of information is obtained assuming that information is lost proportionally to the excess area of the numerical area of influence when compared to the true area of influence. }
\label{fig:L2_rel_error}
\end{figure}

\subsection{Full adjoint-based optimization results}
\label{sec:numerical}

In order to quantify the effects of the aforementioned numerical diffusion on the full adjoint-based optimization,  the adjoint framework in GEOS-Chem was used to produce an optimal reconstruction of the plume's location and magnitude through successive iterations. The best reconstructed plumes produced by the optimization algorithm \citep{bib:Zhu94} after 99 iterations are shown, for both spatial resolutions, in Figures \ref{fig_reonstructions_4x5} and \ref{fig_reonstructions_2x25} for the 20th pressure level, and in Figures \ref{fig_reconstr_4x5_lat_vs_alt} and \ref{fig_reconstr_2x25_lat_vs_alt} for a normalized zonally-integrated side-view. The total number of iterations was considered appropriate based on the limits of computational resources in practical situations. 
For a time simulation of $t_f=3$ hr, the adjoint-based optimization produced excellent reconstructions for both spatial resolutions. Note however, that as the total simulation time increased in our simulations, the quality of the reconstructions decreased rapidly. This is shown in the multiple panels in Figures \ref{fig_reonstructions_4x5}, \ref{fig_reonstructions_2x25}, \ref{fig_reconstr_4x5_lat_vs_alt}, and \ref{fig_reconstr_2x25_lat_vs_alt}. \\

Quantitatively speaking, since the adjoint-based optimization method is formulated in the inner-product space of square integrable functions, it is natural to use the induced $L^2$-norm to evaluate the performance of the method. In Figure \ref{fig:L2_rel_error} we show the relative $L^2$-error between the true initial plume and the reconstructed plume as a function of simulation time for both spatial resolutions. For short simulation times, the reconstructions are excellent. As the simulation time increases, the relative error increases very rapidly. In fact, according to this metric, relative errors above 70$\%$ are incurred in both resolutions in 2-day simulations.

\subsection{Estimation of the loss of information}
\label{sec:loss_of_info}

A simple calculation to estimate of the relative error observed in the adjoint-based source reconstructions is proposed here based on the horizontal deformation of the numerical region of influence, due to undesired numerical diffusive processes, as a function of time.  Estimates of the numerical broadening of the area of influence, for GEOS-Chem, were discussed in section \ref{sec:geo_inf}. In the approach proposed here, it is assumed that the loss of information (specifically our inability to reconstruct the initial condition accurately) increases as the numerically reconstructed area of influence, $A_h(t_0)$, spreads when compared to the original area of the plume $A(t_0)$. Thus, if the numerical area of influence has broadened two-fold ($A_h(t_0)=2A(t_0)$), it is assumed that one can only recover 50\% of the original information (of the location and magnitude of the initial condition) through successive iterations, leading to a 50\% relative error in the adjoint-based reconstruction, calculated as 
\[\text{relative error} (\%)=100 \times (1-A(t_0)/A_h(t_0)).\] 
If the numerical area of influence has tripled ($A_h(t_0)=3A(t_0)$), then one can only recover 33\% of the information, leading to a 66\% relative error, and so on. 
This estimate is plotted in Figure \ref{fig:L2_rel_error}.
The loss of information estimates for both $4^{\circ}\times 5^{\circ}$ and $2^{\circ}\times 2.5^{\circ}$ spatial resolutions are very similar and thus only one representative curve was plotted. Figure \ref{fig:L2_rel_error} illustrates that this simple calculation captures well the salient features of the experimental results obtained in the numerical adjoint-based reconstructions. Note that for simulations of less than 2-days, the adjoint-based iterative process produces better results than the proposed estimate. For longer simulation times the estimate seems appropriate if not slightly optimistic.\\

The quality of the gradient information provided by the adjoint operator is more accurate, as discussed before, for shorter simulations. As a consequence, the BFGS gradient-based minimization routine utilized in GEOS-Chem \citep{bib:Zhu94}, is capable of producing excellent results by decreasing the normalized value of the cost function by 15 orders of magnitude for the $4^{\circ}\times 5^{\circ}$, 3 hr simulations, and nearly by 14 orders of magnitude for the $2^{\circ}\times 2.5^{\circ}$, 3 hr simulations. This is shown in Figures \ref{fig:cost4x5} and \ref{fig:cost2x25}. As the simulation time increases though, the gradient information loses quality rapidly and leads to significantly poorer performances in the minimization processes. Indeed, for 7-day simulations, the minimization only achieves an 8-order of magnitude reduction of the normalized cost function, leading to very poor reconstructed plumes. \\

\begin{figure}[t]
\noindent\includegraphics[width=.5\textwidth]{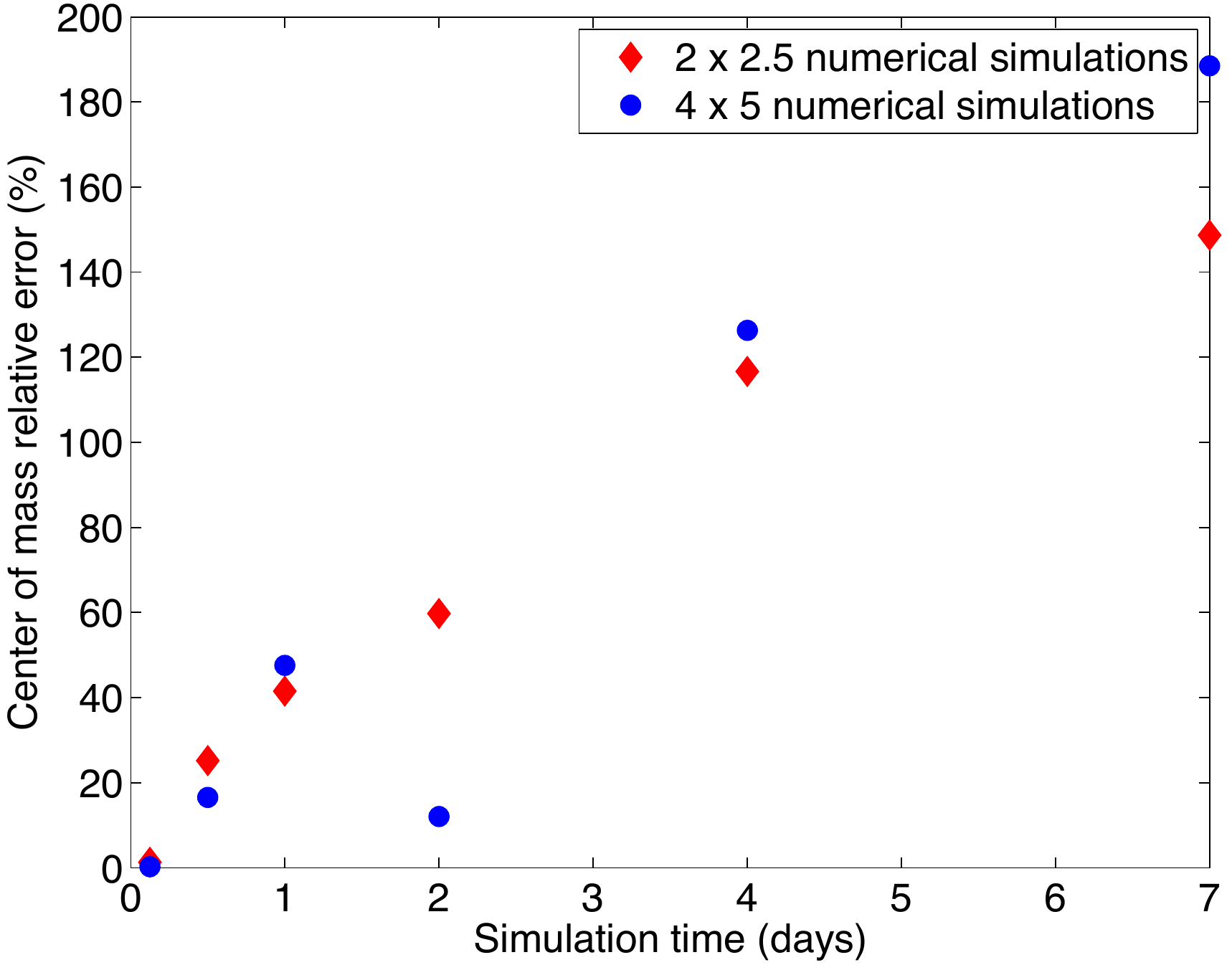}
\caption{Center of mass relative error as a function of simulation time. $100\%$ error corresponds to a reconstructed center of mass at a distance equal to the maximal radius of the plume ($1680$ km) away from the true center of mass.}
\label{fig_center_error}
\end{figure}

\begin{figure}[t]
\noindent\includegraphics[width=.5\textwidth]{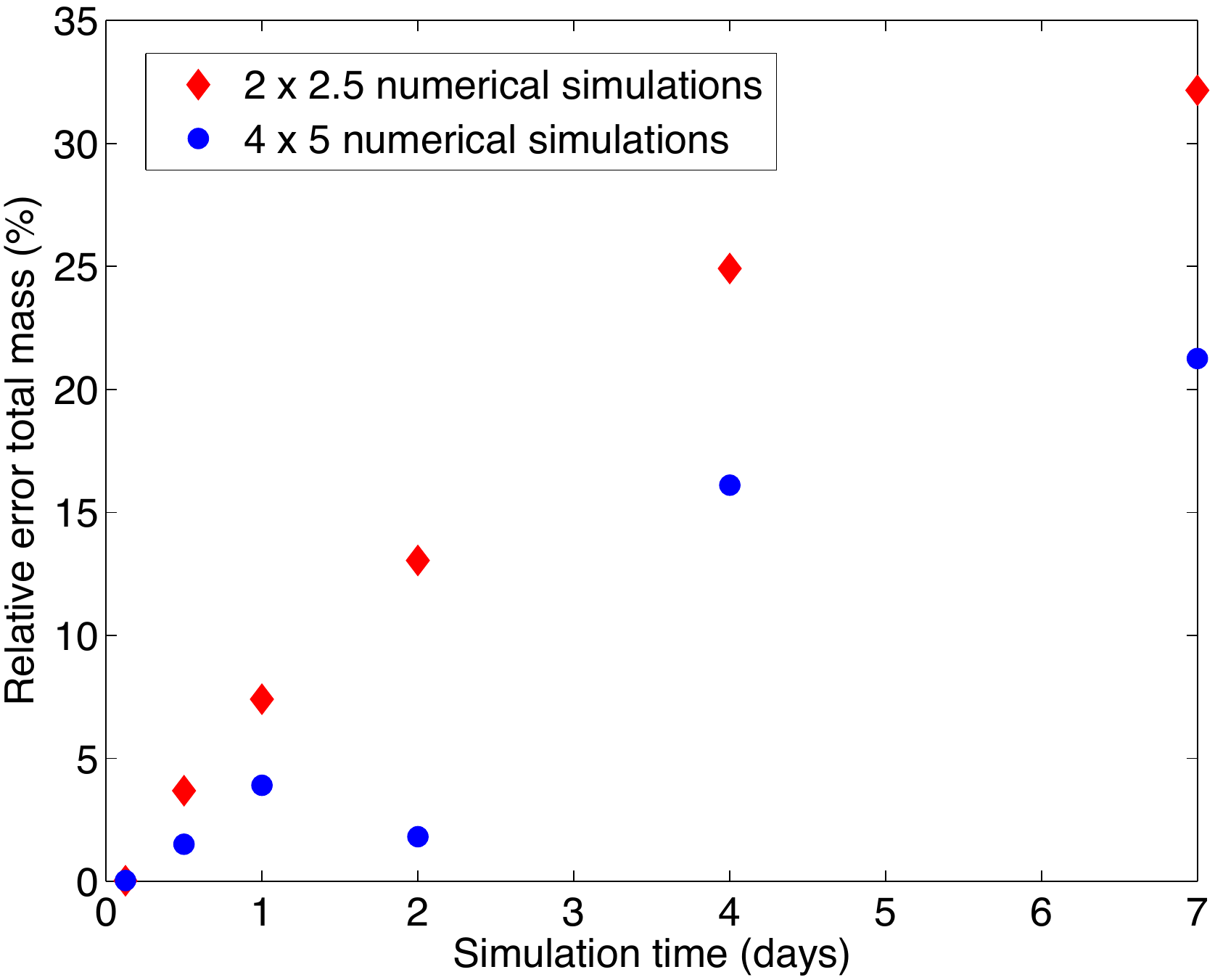}
\caption{Total mass relative error as a function of simulation time.}
\label{fig_total_mass}
\end{figure}

The center of mass and total mass relative errors were calculated, between the reconstructed plume and the true plume, for each spatial resolution and simulation period, in order to provide a more complete characterization of the numerical errors. The results are shown in  Figures \ref{fig_center_error} and \ref{fig_total_mass}. For these calculations, an error in the center of mass equal to the maximal radius of the plume ($1680$ km) was considered a $100\%$ relative error. Errors above $100\%$ are observed in simulations longer that 2 days. The total mass relative error does not reach very high values, for 7-day simulations these are of the order of $20-30\%$, as shown in Figure \ref{fig_total_mass}, due to the conservative properties of the advection scheme. In simulations longer than 2 days, the total plume mass was consistently over-estimated. In all of the numerical experiments unphysical negative concentrations in the reconstructed plume were not observed. Figures \ref{fig_reonstructions_4x5}, \ref{fig_reonstructions_2x25}, \ref{fig_reconstr_4x5_lat_vs_alt}, and \ref{fig_reconstr_2x25_lat_vs_alt} show that the reconstructed plume is deformed both in horizontal and the vertical components. This is due to vertical atmospheric convective processes represented as sub-grid parametrization in GEOS-Chem. \\

\subsection{Relevance in continuous-in-time source reconstructions}
In the context of reconstruction of continuous-in-time sources or top-down flux source attribution problems, one way to interpret the implications of the local-in-time findings presented here, is to calculate the product of the time scale, associated with a given relative error upper limit, and the mean flow velocity, to find the distance between the continuous sources and observation sites, beyond which measurements may not successfully be used for the reconstruction of the magnitude and location of the source in Eulerian CTMs. This reconstruction will be compromised since the signal to noise ratio will be smaller and smaller in subsequently further observation sites. In the experiments presented here, assuming a mean atmospheric flow of 10 $m$ $s^{-1}$, and a 2 day time scale, the resulting length scale for which measurements may not be successfully used to reconstruct a source is about 1700 $km$, equivalent to 3 grid boxes in a $4^{\circ}\times 5^{\circ}$ simulation. This means that if a continuous source is to be successfully reconstructed using the adjoint-based optimization method in current global CTMs, then a continuous observation system, closer than 3 grid cells (in a $4^{\circ}\times 5^{\circ}$ simulation), should be in place down wind from the source. Moreover, the aforementioned length scale naturally suggests a spatial density over the globe below which sources may not be accurately reconstructed in global Eulerian CTMs in practical spatial resolutions. Depending on the tracer, this spatial density of observation sites may or may not be realistic with current Earth observing systems. \\

\section{Conclusions}

It is shown that the ability to reconstruct the location and magnitude of an accident-type, instantaneous inert pollution plume, using the adjoint-based optimization method in a global Eulerian CTM, decreases rapidly as the the time between release and observations increases. The obtained numerical results suggest a time scale of 2 days after which significant numerical errors ($> 70\%$) compromise the ability of the adjoint-based optimization method to track information accurately back to the pollution source. This time scale is shorter than the time scales of inter-continental transport for example, for which it would be desirable to have a reliable methodology to reconstruct sources. \\

It is shown that the quality of the gradient utilized to minimize the misfit between observations and simulation results, used by the adjoint-based optimization method, decays fast as a function of simulation time. As a consequence of this fact, a simple way to quantitatively characterize the loss of information in the adjoint-based optimization is proposed based on the --undesired-- deformation of the numerical area of influence caused by the effective diffusion in the CTM (which is much larger than expected by the order of the advection scheme due to the chaotic nature of atmospheric flow). This characterization of the loss of information is shown to successfully describe the behavior of the numerical relative error observed between the true solution and the numerical reconstructions, as a function of simulation time.\\

The approach presented in this study exhibits a structural (algorithmical) limitation also present in data assimilation systems utilizing the adjoint to obtain the gradient for successive optimization iterations in source attributions problems. Since identical twin experiments were utilized to evaluate the performance of the adjoint-based optimization methodology, and measurement errors and model errors were ignored, one should expect the results obtained here to be optimistic in the context of real-life source attribution problems.  \\

The use of regularization techniques as well as a dense observation network near the source(s) to be reconstructed may help counter balance this structural numerical limitation of the adjoint method, as it is currently implemented in global CTMs. Inversion approaches not relying on gradient information should be more intensely studied in this context \citep{bib:Kai04}. \\\\




\begin{figure}[t]
\noindent\includegraphics[width=.5\textwidth]{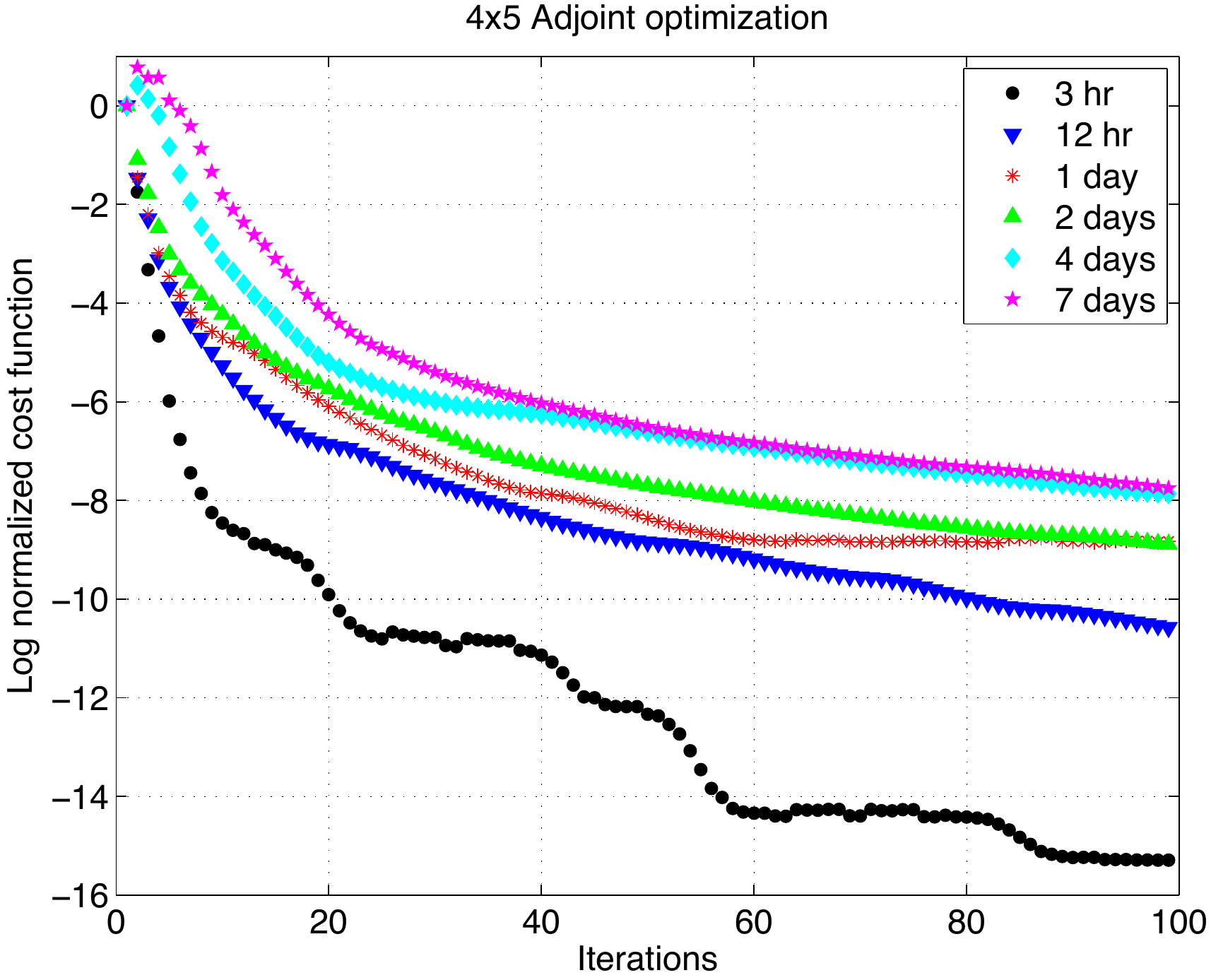}
\caption{Behavior of the cost function as a function of the number of iterations for multiple simulations times for $4^{\circ} \times 5^{\circ}$ spatial resolution. Logarithm of the normalized cost function vs number of iterations.}
\label{fig:cost4x5}
\end{figure}

\begin{figure}[t]
\noindent\includegraphics[width=.5\textwidth]{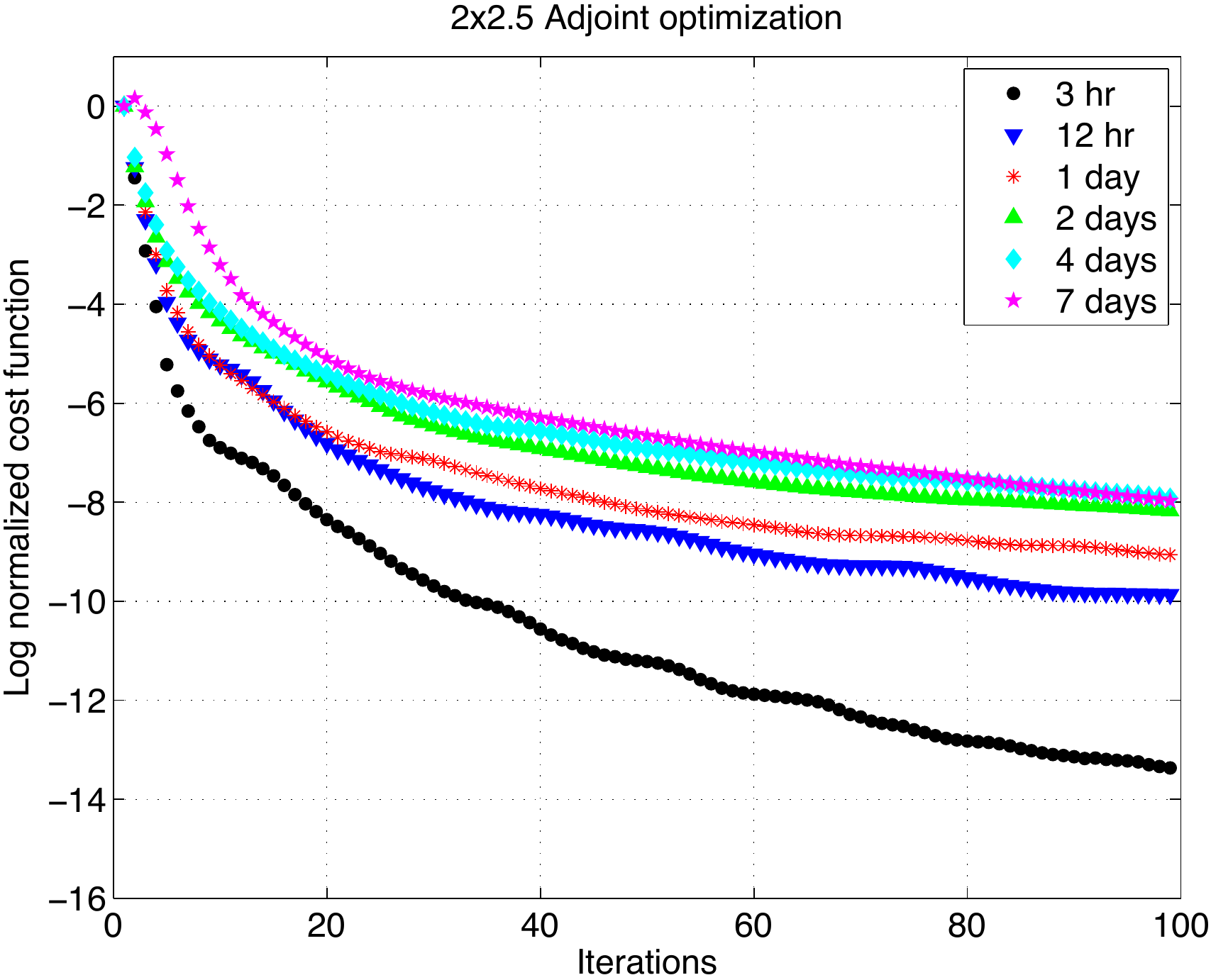}
\caption{Behavior of the cost function as a function of the number of iterations for multiple simulations times for $2^{\circ} \times 2.5^{\circ}$ spatial resolution. Logarithm of the normalized cost function vs number of iterations.}
\label{fig:cost2x25}
\end{figure}

\begin{acknowledgment} 

\end{acknowledgment}

\ifthenelse{\boolean{dc}}
{}
{\clearpage}
%
%
%
%
%
%

\ifthenelse{\boolean{dc}}
{}
{\clearpage}


%
\end{document}